\documentclass{article}
\usepackage{epic,latexsym,amssymb,soul,xcolor}
\usepackage{color}
\usepackage{tikz}
\usepackage{amsfonts,epsf,amsmath}

\newtheorem{thm}{Theorem}
\newtheorem{lem}{Lemma}

\newtheorem{claim}{Claim}

\newtheorem{conj}{Conjecture}

\newcommand{\qed}{$\Box$}

\newcommand{\proof}{\noindent\textbf{Proof. }}

\let\oldenumerate\enumerate
\renewcommand{\enumerate}{
  \oldenumerate
  \setlength{\itemsep}{0pt}
  \setlength{\parskip}{0pt}
  \setlength{\parsep}{0pt}
}

\begin{document}

\title{Disjoint Cross Intersecting Families}
\markright{Disjoint Cross Intersecting Families}

\author{Nuttanon Songsuwan, Supida Sengsamak, Nutchapol Jeerawattana,\\
 Thiradet Jiarasuksakun \, and \, Pawaton Kaemawichanurat$^\dag$
\\ \\
Department of Mathematics, Faculty of Science,\\
King Mongkut's University of Technology Thonburi, \\
$^{\dag}$Discrete Mathematics, Algebra, Probability, Statistics, \\
Dynamical Systems and Applications (DAPSDAP), \\
Bangkok, Thailand \\
\small \tt Email: nuttanon.19701@gmail.com, nutchapolj@gmail.com, supida.gift@gmail.com, \\
thiradet.jia@kmutt.ac.th, pawaton.kae@kmutt.ac.th}

\date{}
\maketitle

\begin{abstract}
\noindent For positive integers $n$ and $r$ such that $r \leq \lfloor n/2\rfloor$, let $X$ be a set of $n$ elements and let $\binom{X}{r}$ be the family of all $r$-subsets of $X$. Two sub-families $\mathcal{A}$ and $\mathcal{B}$ of $\binom{X}{r}$ are called cross intersecting if $A \cap B \neq \emptyset$ for all $A \in \mathcal{A}$ and $B \in \mathcal{B}$. One of main tools in the study of extremal set theory, and cross intersecting families in particular, is compression operation. In this paper, we give an example of cross intersecting families $\mathcal{A}$ and $\mathcal{B}$ that the compression operation is not applicable when $\mathcal{A}$ and $\mathcal{B}$ are disjoint. We develop new technique to prove that, for disjoint cross intersecting families $\mathcal{A}$ and $\mathcal{B}$ of $\binom{X}{r}$,
\begin{displaymath}
|\mathcal{A}| + |\mathcal{B}| \leq \binom{n}{r} - \binom{l}{p}
\end{displaymath}
where $n = 2r + l$ and $p = min\{r, \lceil \frac{l}{2}\rceil\}$. This bound is asymptotically sharp.
\end{abstract}


%\vskip 1 cm
%\newpage
\section{Introduction}
For a natural number $n$, we let $[n] = \{ 1,2,\dots,n \}$. We let $\binom{[n]}{r} = \{ A \subseteq [n] : |A| = r \}$. For a subset set $T$ of $[n]$ where $|T| = t$, a subfamily $\mathcal{S}$ of $\binom{[n]}{r}$ is said to be \textit{t-star} if $T \subseteq S$ for all $S \in \mathcal{S}$. In particular, if $|T|=1$, then $\mathcal{S}$ is simply called \textit{a star}. A subfamily $\mathcal{A}$ of $\binom{[n]}{r}$ is said to be \emph{intersecting} if $A \cap B \neq \emptyset$ for all $A, B \in \mathcal{A}$. Let $\mathcal{A}$ and $\mathcal{B}$ be subfamilies of $\binom{[n]}{r}$. Then $\mathcal{A}$ and $\mathcal{B}$ are called $t$-\textit{cross intersecting families} if $|A \cap B| \geq t$ for any $A \in \mathcal{A}, B \in \mathcal{B}$. In particular, $\mathcal{A}$ and $\mathcal{B}$ are called \emph{cross intersecting families} when $t = 1$. Further, cross intersecting families $\mathcal{A}$ and $\mathcal{B}$ are said to be disjoint if $\mathcal{A} \cap \mathcal{B} = \emptyset$.
\vskip 5 pt

%\indent A graph $G$ consists of the vertex set $V(G)$ and the edge set $E(G) \subseteq \binom{V(G)}{2}$. For a subfamily $\mathcal{A}$ of $\binom{V}{r}$, the \emph{intersecting} graph $G_{\mathcal{A}}$ is the graph with $V(G_{\mathcal{A}}) = \mathcal{A}$ and vertices $X, Y \in V(G_{\mathcal{A}})$ are joined by an edge if and only if $X \cap Y \neq \emptyset$. An \emph{independent set} $S$ is a set of vertices which are pairwise non-adjacent. A \emph{complete} $k$\emph{-partite graph} is the graph $G$ whose $V(G)$ is partitioned into $k$ independent sets such that each vertex in each set is adjacent to all vertices in the other sets.

\indent One of the famous problems in extremal set theory is to determine how large a subfamily $\mathcal{A}$ of $\binom{[n]}{r}$ can be when the family $\mathcal{A}$ is intersecting. When $\frac{n}{2} < r \le n$, the problem become trivial since for any $A,B \in \binom{[n]}{r}$, $n \ge |A \cup B| = |A| +|B| - |A \cap B| = 2r - |A \cap B| > n - |A \cap B|$ which implies $|A \cap B| > 0$. Hence, $|\mathcal{A}| \leq \binom{n}{r}$. When $1 \le r \le \frac{n}{2}$, the solution was classical result, called \emph{Erd\H os-Ko-Rado Theorem}~\cite{Erdos1961}.

\begin{thm}[Erd\H os-Ko-Rado Theorem]~\cite{Erdos1961}\label{erdos}
If $\mathcal{A}$ is an intersecting family of $\binom{[n]}{r}$, then
\begin{displaymath}
|\mathcal{A}| \le \binom{n-1}{r-1}
\end{displaymath}
where $1 \le r \le \frac{n}{2}$. The equality holds if $\mathcal{A}$ is a star of $\binom{[n]}{r}$.
\end{thm}

\indent Pyber~\cite{Pyber1986} extended this study to find the upper bound of product size of cross intersecting families and this has been extensively studied see~\cite{Borg2014,Borg2015,Borg2017,Matsumoto1989} for example. In 2010, Tokushige~\cite{Norihide2010} generalized Pyber's result to $t$-cross intersecting families.

\begin{thm}\label{PN}\cite{Pyber1986,Norihide2010}
If $\mathcal{A}$ and $\mathcal{B}$ are $t$-cross intersecting families of $\binom{[n]}{r}$ for $t \geq 1$, then
\begin{displaymath}
|\mathcal{A}||\mathcal{B}| \leq \binom{n - 1}{r - 1}^{2}
\end{displaymath}
where $1 \le r \le \frac{n}{2}$. The equality holds if $\mathcal{A} = \mathcal{B}$ is a $t$-star of $\binom{[n]}{r}$.
\end{thm}

\indent For the study when cross intersecting families $\mathcal{A}$ and $\mathcal{B}$ are disjoint, Kupaskii~\cite{Kupavskii} conjectured that :

\begin{conj}\cite{Kupavskii}\label{Ku}
If $\mathcal{A}$ and $\mathcal{B}$ are disjoint cross intersecting families of $\binom{[n]}{r}$, then
\begin{displaymath}
\min\{|\mathcal{A}|, |\mathcal{B}|\} \leq \frac{1}{2}\binom{n - 1}{r - 1}.
\end{displaymath}
\end{conj}

\noindent Very recently, Conjecture \ref{Ku} was proved by Huang~\cite{Huang2019} when $r \geq 2$ and $n \geq 2r^{2}$ and was disproved when $n < r^{2}$. Frankl and Kupavskii~\cite{Frankl2020} independently proved that Conjecture \ref{Ku} is true when $n \geq cr^{2} + (4 - 2c)r$ and disproved this conjecture when $n \leq cr^{2} - 2cr + 1$ where $c = \log_{2}e$.
\vskip 5 pt

\indent In the study of extremal set theory, there is a useful technique called \emph{compression operation}, or \emph{shifting}.
%The technique is used to organized almost all set members of the family to have less number of common elements but still preserve the size of family as well as intersecting property, and then, we apply induction.
The technique was introduced in~\cite{Erdos1961} and is detailed as follows. Let $\mathcal{A}$ be an intersecting family and let $A \in \mathcal{A}$. For $i, j \in [n]$, we let

\[ \delta_{i, j}(A)  = \left\{\begin{array}{lll}
                       (A \setminus \{j\}) \cup \{i\} & if~~j \in A~and~i \notin A,\\
                       A                              & otherwise.\\
               \end{array}\right. \]

\noindent Then we let,
\vskip 5 pt

\indent $\Delta_{i, j}(\mathcal{A}) = \{A \in \mathcal{A} : \delta_{i, j}(A) \in \mathcal{A}\} \cup \{\delta_{i, j}(A) :  A \in \mathcal{A}, \delta_{i, j}(A) \notin \mathcal{A}\}$.
\vskip 5 pt

\noindent We apply this operation for all $i < j$ repeatedly until we have the family $\hat{\mathcal{A}}$ such that $\delta_{i, j}(\hat{\mathcal{A}}) = \hat{\mathcal{A}}$ for all $i < j$. By this operation, it can be proved that $(1)$ $|\hat{\mathcal{A}}| = |\mathcal{A}|$, $(2)$ $\hat{\mathcal{A}}$ is an intersecting family and $(3)$ for each $A, B \in \hat{\mathcal{A}}$ we have that $A \cap B \cap [n - 1] \neq \emptyset$. We see that, by $(3)$, we can apply induction on the sets of $\hat{\mathcal{A}}$.
\vskip 5 pt

\indent It has been convinced, by Borg~\cite{Borg2015} and Togushige~\cite{Norihide2010} for example, that for cross intersecting families $\mathcal{A}$ and $\mathcal{B}$ of $\binom{[n]}{r}$, the families

\indent $\Delta_{i, j}(\mathcal{A}) = \{A \in \mathcal{A} : \delta_{i, j}(A) \in \mathcal{A}\} \cup \{\delta_{i, j}(A) : A \in \mathcal{A}, \delta_{i, j}(A) \notin \mathcal{A}\}$,
\vskip 5 pt

\indent $\Delta_{i, j}(\mathcal{B}) = \{B \in \mathcal{B} : \delta_{i, j}(B) \in \mathcal{B}\} \cup \{\delta_{i, j}(B) : B \in \mathcal{B}, \delta_{i, j}(B) \notin \mathcal{B}\}$.
\vskip 5 pt

\noindent are still cross intersecting. However, this is not always true when $\mathcal{A}$ and $\mathcal{B}$ are disjoint.
\vskip 10 pt

\section{Main results}
In this section, we first give an example of disjoint cross intersecting families that the compression operation is not applicable. Let $\mathcal{A}$ and $\mathcal{B}$ be disjoint cross intersecting families of $\binom{[n]}{r}$ when $2r \leq n$. Suppose, up to isomorphism, that $A = \{1, r + 2, ..., 2r\} \in \mathcal{A}$ and $B = \{2, r + 2, ..., 2r\} \in \mathcal{B}$. If $i = 1$ and $j = 2$, then we have that $\delta_{1, 2}(A) = A$ and $\delta_{1, 2}(B) = A$. Since $A \in \mathcal{A}$ and $A \notin \mathcal{B}$, it follows that $\delta_{1, 2}(B) \in \Delta_{1, 2}(\mathcal{A}) \cap \Delta_{1, 2}(\mathcal{B})$ violating the disjoint property. Even though, we further modify the definition of $\Delta_{i, j}(\mathcal{A})$ and $\Delta_{i, j}(\mathcal{A})$ to $\delta_{i, j}(A), \delta_{i, j}(B) \in \mathcal{A} \cup \mathcal{B}$ and $\delta_{i, j}(A), \delta_{i, j}(B) \notin \mathcal{A} \cup \mathcal{B}$ in order to preserve the disjoint property between $\Delta_{i, j}(\mathcal{A})$ and $\Delta_{i, j}(\mathcal{B})$, the intersecting property collapses. We let
\vskip 5 pt

\indent $\Delta_{i, j}(\mathcal{A}) = \{A \in \mathcal{A} : \delta_{i, j}(A) \in \mathcal{A} \cup \mathcal{B}\} \cup \{\delta_{i, j}(A) : A \in \mathcal{A}, \delta_{i, j}(A) \notin \mathcal{A} \cup \mathcal{B}\}$,
\vskip 5 pt

\indent $\Delta_{i, j}(\mathcal{B}) = \{B \in \mathcal{B} : \delta_{i, j}(B) \in \mathcal{A} \cup \mathcal{B}\} \cup \{\delta_{i, j}(B) : B \in \mathcal{B}, \delta_{i, j}(B) \notin \mathcal{A} \cup \mathcal{B}\}$.
\vskip 5 pt

\noindent Recall that  $A = \{1, r + 2, ..., 2r\} \in \mathcal{A}$ and $B = \{2, r + 2, ..., 2r\} \in \mathcal{B}$. Moreover, we suppose that $\{1, 3, ..., r + 1\} \notin \mathcal{A} \cup \mathcal{B}$ and $C = \{2, 3, ..., r + 1\} \in \mathcal{A}$. If $i = 1$ and $j = 2$, then
\vskip 5 pt

\indent $\delta_{1, 2}(C) = \{1, 3, ..., r + 1\} \notin \mathcal{A} \cup \mathcal{B}$ and $\delta_{1, 2}(B) = A \in \mathcal{A} \cup \mathcal{B}$.
\vskip 5 pt

\noindent Therefore, $\delta_{1, 2}(C) \in \Delta_{1, 2}(\mathcal{A})$ and $B \in \Delta_{1, 2}(\mathcal{B})$, but $\delta_{1, 2}(C) \cap B = \emptyset$. So the compression operation does not preserve intersecting property.
\vskip 5 pt

\indent Hence, in this paper, we develop new technique~(see Claim 1) to prove that

\begin{thm}\label{thm main}
For disjoint cross intersecting families $\mathcal{A}$ and $\mathcal{B}$ of $\binom{[n]}{r}$,
\begin{displaymath}
|\mathcal{A}| + |\mathcal{B}| \leq \binom{n}{r} - \binom{l}{p}
\end{displaymath}
when $n = 2r + l$ and $p = min\{r, \lceil \frac{l}{2}\rceil\}$.
\end{thm}

\noindent The proof of Theorem \ref{thm main} is given in the next section. For sharpness of the bound in this theorem, when $r$ is small compare to $n$, we fix $r$ and let $n$ be sufficiently large. We have that $p = min\{r, \lceil\frac{l}{2}\rceil\} = r$ and $l = n - 2r$. Thus, the upper bound of Theorem \ref{thm main} is $\binom{n}{r} - \binom{l}{p} = \binom{n}{r} - \binom{n - 2r}{r} = \Theta(n^{r - 1})$. Let $\mathcal{S}$ be a star and let $\mathcal{A}_{1}, \mathcal{B}_{1}$ be a partition of $\mathcal{S}$. That is, $\mathcal{A}_{1}$ and $\mathcal{B}_{1}$ are disjoint cross intersecting families and $|\mathcal{A}_{1}| + |\mathcal{B}_{1}| = |\mathcal{S}| = \binom{n - 1}{r - 1} = \Theta(n^{r - 1})$. It is worth noting that, the upper bound of $|\mathcal{A}||\mathcal{B}|$ in this case is maximized at $\frac{(\binom{n}{r} - \binom{l}{p})^{2}}{4}$ which is $\Theta(n^{2r - 2})$ the same as that of Theorem \ref{PN} when $t = 1$.
\vskip 5 pt

\indent When $r$ is large compare to $n$, we fix $l$ and let $n$ be sufficiently large, we have that $p = min\{r, \lceil\frac{l}{2}\rceil\} = r$. Thus, the upper bound of Theorem \ref{thm main} is $\binom{n}{r} - \binom{l}{p} = \Theta(n^{r})$. For the construction, we let

\begin{displaymath}
\mathcal{X} = \{ A \in \binom{[2r]}{r} : 1 \in A \textnormal{ and } 2 \notin A \} \text{ and let } \mathcal{X}^c = \{[2r] \setminus A : A \in \mathcal{X} \}.
\end{displaymath}

\noindent Clearly, for any $A^{c} \in \mathcal{X}^{c}$, we have that $|A^{c}| = r, 1 \notin A^{c}$ and $2 \in A^{c}$. Thus, $\mathcal{X} \cap \mathcal{X}^c = \emptyset$ and $|\mathcal{X}^c|=|\mathcal{X}| = \binom{2r-2}{r-1}$. Then, we let
\vskip 5 pt

\begin{displaymath}
\mathcal{A}_{1} = \mathcal{X} \cup \mathcal{X}^c.
\end{displaymath}

\noindent Thus $1 \in A$ and $2 \notin A$ or $1 \notin A$ and $2 \in A$ for any $A \in \mathcal{A}_{1}$. Moreover, we let
\vskip 5 pt

\begin{displaymath}
\tilde{\mathcal{B}}_{1} = \{B \in \binom{[2r]}{r} : 1, 2 \notin B\} \text{ and } \tilde{\mathcal{B}}_{2} =  \{B \in \binom{[n]}{r} : 1, 2 \in B\}.
\end{displaymath}

\noindent Then, we let
\vskip 5 pt

\begin{displaymath}
\mathcal{B}_{1} = \tilde{\mathcal{B}}_{1} \cup \tilde{\mathcal{B}}_{2}.
\end{displaymath}

\noindent Clearly, $\mathcal{A}_{1}$ and $\mathcal{B}_{2}$ are disjoint. Let $A \in \mathcal{A}_{1}$ and $B \in \mathcal{B}_{2}$. If $B \in \tilde{\mathcal{B}}_{1}$, then $B \neq ([2r] \setminus A)$. Thus, $B \cap A \neq \emptyset$. If $B \in \tilde{\mathcal{B}}_{2}$, then $1 \in A \cap B$ or $2 \in A \cap B$. Therefore, $\mathcal{A}_{1}$ and $\mathcal{B}_{1}$ are disjoint cross intersecting family. Clearly,

\begin{displaymath}
|\mathcal{A}_{1}| = |\mathcal{X}| + |\mathcal{X}^c| = 2\binom{2r-2}{r-1} \text{ and } |\mathcal{B}_{1}| = \binom{2r - 2}{r} + \binom{n - 2}{r - 2}.
\end{displaymath}

\noindent Because $2r = n - l$, it follows that

\begin{align}
|\mathcal{A}_{1}| + |\mathcal{B}_{1}| &= 2\binom{2r-2}{r-1} + \binom{2r - 2}{r} + \binom{n - 2}{r - 2}\notag\\
                                      &= \binom{2r-2}{r-1} + \Big(\binom{2r-2}{r-1} + \binom{2r - 2}{r}\Big) + \binom{n - 2}{r - 2}\notag\\
                                      &= \binom{2r-2}{r-1} + \binom{2r - 1}{r} + \binom{n - 2}{r - 2}\notag\\
                                      &= \binom{n - l - 1}{r} + \binom{n - l - 2}{r-1} + \binom{n - 2}{r - 2} = \Theta(n^{r}).\notag
\end{align}
\vskip 5 pt

\indent Finally, we observe that when $n = 2r$ (that is when $l = 0$), we have $|\mathcal{A}| + |\mathcal{B}| \leq \binom{n}{r}$. We will characterize disjoint cross intersecting families $\mathcal{A}_{1}$ and $\mathcal{B}_{1}$ of $\binom{[2r]}{r}$ such that $|\mathcal{A}_{1}| + |\mathcal{B}_{1}| = \binom{2r}{r}$. We need the following lemma in characterization.
\vskip 5 pt

\begin{lem}\label{complete-k}
Let $G_{1}$ be the graph having vertex set $\binom{[2r]}{r}$ and edge set $E = \{XY : X, Y \in \binom{[2r]}{r}$ and $X \cap Y \ne \emptyset\}$. Then $G_{1}$ is a complete $k$-partite graph where $k=\frac{1}{2}\binom{2r}{r} = \binom{2r - 1}{r}$.
\end{lem}
\proof Let $k = \binom{2r-1}{r - 1} = \binom{2r-1}{r}$ and $\{X_1,X_2, \dots, X_k\} = \{X:X \in \binom{[2r-1]}{r} \}$. For $1 \le i \le k$, let $Y_i = [2r] \backslash X_i$. Since $2r-1 < 2r$, it follows that $X_i \cap X_j \ne \emptyset$, for $1 \le i \ne j \le k$. Moreover $2r \in Y_i$. So $Y_i \cap Y_j \ne \emptyset$. By the definition of $Y_i$, $X_i \cap Y_i = \emptyset$.
\vskip 5 pt

\indent Let $V_i = \{X_i, Y_i\}$. We will show that $G_{1}$ is a complete $k$-partite graph whose partite sets are $V_1,V_2, \dots,$ $V_k$. Since $X_i \cap Y_i = \emptyset$, there is no edge in $V_i$. Because $X_i \cap X_j \ne \emptyset$, $X_i$ is adjacent to $X_j$. Since $X_i \ne X_j$, there exists $a \in X_i \backslash X_j$. So $a \in [2r] \backslash X_j = Y_j$. Therefore, $X_i \cap Y_j \ne \emptyset$. Hence, $X_i$ is adjacent to $Y_j$. This completes the proof.
%\qed
\qed
\vskip 5 pt

\noindent By Lemma \ref{complete-k}, we let $V_{1}, V_{2}, ..., V_{k}$ be the partite sets of $G_{1}$. For some $0 \leq t \leq k$, we let $\mathcal{A}_{1}$ be the union of $V_{i_{1}}, V_{i_{2}}, ..., V_{i_{t}}$ where $i_{1}, i_{2}, ..., i_{t}$ is a subsequence of $1, 2, .., k$. Then, we let $\mathcal{B}_{1}$ be $\cup^{k}_{i = 1}V_{i} \setminus \cup^{t}_{j = 1}V_{i_{j}}$. We characterize that
\vskip 5 pt

\begin{thm}\label{2kCk}
Let $\mathcal{A},\mathcal{B}$ be disjoint cross intersecting families of $\binom{[2r]}{r}$. Then $|\mathcal{A}| + |\mathcal{B}| \leq \binom{2r}{r}$ and equality holds if and only if the families $\mathcal{A}, \mathcal{B}$ are $\mathcal{A}_{1}, \mathcal{B}_{1}$.
\end{thm}
\proof By Theorem \ref{thm main}, we establish the upper bound. Hence, we assume that $|\mathcal{A}| + |\mathcal{B}| = \binom{2r}{r}$. Because $\mathcal{A}$ and $\mathcal{B}$ are disjoint, it follows that

\begin{displaymath}
\mathcal{A} \cup \mathcal{B} = \binom{[2r]}{r}.
\end{displaymath}

\noindent Thus $\mathcal{A} \cup \mathcal{B} = \{X_{1}, ..., X_{k}, Y_{1}, ..., Y_{k}\}$, where $X_{i}$ and $Y_{i}$ are defined in Lemma \ref{complete-k} and $k = \binom{2r-1}{r-1}$. By intersecting property, we have that $X_{i} \in \mathcal{A}$ if and only if $Y_{i} \in \mathcal{A}$. For $1 \leq i \leq k$, let $V_{i} = \{X_{i}, Y_{i}\}$. So there exists a subsequence $V_{i_{1}}, V_{i_{2}}, ..., V_{i_{t}}$ such that $\mathcal{A} = \cup^{i}_{j = 1}V_{i_{j}}$ and $\mathcal{B} = \cup^{k}_{i = 1}V_{i} \setminus \cup^{i}_{j = 1}V_{i_{j}}$. Therefore, the families $\mathcal{A}$ and $\mathcal{B}$ are $\mathcal{A}_{1}$ and $\mathcal{B}_{1}$, respectively. This completes the proof.
\qed

\section{Proof of Theorem \ref{thm main}}
By Theorem \ref{2kCk}, we consider when $\binom{l}{p} \neq 0$. We will show that $|\mathcal{A}| + |\mathcal{B}| < \binom{n}{r} - \binom{l}{p}$. Suppose to the contrary that
\begin{align}\label{eq 1-3}
|\mathcal{A}| + |\mathcal{B}| \geq \binom{n}{r} - \binom{l}{p}.
\end{align}
\noindent Thus, there are at most $\binom{l}{p}$ sets which are in $\binom{[n]}{r} \setminus (\mathcal{A} \cup \mathcal{B})$. We may let $\mathcal{C} = \binom{[n]}{r} \setminus (\mathcal{A} \cup \mathcal{B})$. Hence, $|\mathcal{C}| \leq \binom{l}{p}$ because $\mathcal{A}$ and $\mathcal{B}$ are disjoint. The following claim is the new technique that we develop and is proved under the assumption (\ref{eq 1-3}).
\vskip 5 pt

\begin{claim}\label{absorb}
If $n = 2r + l$ and $A, B \in \binom{[n]}{r} \setminus \mathcal{C}$ with $A \neq B$, then there exist $S_{0}, S_{1}, ..., S_{f}(f \geq 1)$ such that $A = S_{0}, B = S_{f}$, and for $k = 0, ..., f - 1, S_{k} \cap S_{k + 1} = \emptyset$ and $S_{k}, S_{k + 1} \notin \mathcal{C}$.
\end{claim}
\proof If $A \cap B = \emptyset$, then $l = 0$ and $A = S_{0}$, $B = S_{1}$. Hence, we may assume that $|A \cap B| = t \geq 1$. We may let
\vskip 5 pt

\indent $A = \{x_{1}, ..., x_{t}, z_{1}, ..., z_{r - t}\}$ and $B = \{x_{1}, ..., x_{t}, y_{1}, ..., y_{r - t}\}$
\vskip 5 pt

\noindent where $x_{1}, ..., x_{t}, z_{1}, ..., z_{r - t}$ and $y_{1}, ..., y_{r - t}$ are $2r - t$ different elements of $[n]$. We may distinguish $2$ cases.
\vskip 5 pt

\noindent \textbf{Case 1 :} $\lceil \frac{l}{2}\rceil \geq r$\\
\indent Thus $p = r$. We may let $S_{0} = A$ and $S_{2} = B$. Hence, $|S_{0} \cup S_{2}| = |S_{0}| + |S_{2}| - |S_{0} \cap S_{2}| = 2r - t$. We choose $r$ elements from $[n] \setminus (S_{0} \cup S_{2})$ to be the set $S_{1}$. Clearly, $|[n] \setminus (S_{0} \cup S_{2})| = 2r + l - (2r - t) = l + t > l$. So, there are $\binom{l + t}{r} > \binom{l}{r}$ possible sets. Because $|\mathcal{C}| \leq \binom{l}{r}$, there exists a set $\{a_{1}, ..., a_{r}\} \in \binom{[n] \setminus (S_{0} \cup S_{2})}{r} \setminus \mathcal{C}$. We let
\vskip 5 pt

\indent $S_{1} = \{a_{1}, ..., a_{r}\}$.
\vskip 5 pt

\noindent By this choice, $S_{1} \cap S_{0} = \emptyset$ and $S_{1} \cap S_{2} = \emptyset$. This proves Case 1.
\vskip 5 pt

\noindent \textbf{Case 2 :} $\lceil \frac{l}{2}\rceil < r$\\
\indent So $p = \lceil \frac{l}{2}\rceil$. Recall that $A = \{x_{1}, ..., x_{t}, z_{1}, ..., z_{r - t}\}$ and $B = \{x_{1}, ..., x_{t}, y_{1}, ..., y_{r - t}\}$. For some integers $0 \leq m$ and $0 \leq q < p$, we let $t = mp + q$. We may partition sets $A$ and $B$ as follows. For $1 \leq i \leq m$, we let
\vskip 5 pt

\indent $A_{i} = \{x_{(i - 1)p + 1}, ...,  x_{ip}\}$. Moreover, we let
\begin{displaymath}
A_{m + 1} =
\begin{cases}
\{x_{mp + 1}, ..., x_{t}, y_{1}, ..., y_{r - t}\} & \textnormal{if } r - t < p - q ,\\
\{x_{mp + 1}, ..., x_{t}, y_{1}, ..., y_{p - q}\} &\textnormal{if } r - t \geq p - q.
\end{cases}
\end{displaymath}
\noindent Clearly, $|A_{1}| = \cdots |A_{m}| = p$ and $|A_{m + 1}| \leq p$. We now consider the set $B$. When $m = 0$, we have that $t = q$. We let $B_{1} = \{x_{1},  ..., x_{q}, z_{1}, ..., ,z_{p - q}\}$. When $m \geq 1$, we let
\vskip 5 pt

\indent $B_{i} = \{x_{(m - i)p + q + 1}, ..., x_{(m - i + 1)p + q}\}$
\vskip 5 pt

\noindent for $1 \leq i \leq m$. Moreover, we let
\begin{displaymath}
B_{m + 1} =
\begin{cases}
\{x_{1}, ..., x_{q}, z_{1}, ..., z_{r - t}\} & \textnormal{if } r - t < p - q ,\\
\{x_{1}, ..., x_{q}, z_{1}, ..., z_{p - q}\} &\textnormal{if } r - t \geq p - q.
\end{cases}
\end{displaymath}

\noindent It is worth noting that $|B_{1}| = \cdots = |B_{m}| = p$ and $|B_{m + 1}| \leq p$.
\vskip 5 pt

\noindent Remind that our goal is to construct sets $S_{0}, ..., S_{f}$. However, these sets will be obtained by adding some elements to sets $C_{1}, ..., C_{2(m + 1) + 1}$ as defined in the following. For $1 \leq i \leq m + 1$, we let
\vskip 5 pt

\indent $C_{2i} = A \setminus \cup^{i}_{j = 1}A_{j}$. Moreover, for $0 \leq i \leq m$, we let
\vskip 5 pt

\indent $C_{2i + 1} = B \setminus \cup^{m - i + 1}_{j = 1}B_{j}$ and $C_{2(m + 1) + 1} = B$.
\vskip 5 pt

\indent We see that $C_{2m} \subseteq C_{2m - 2} \subseteq \cdots \subseteq C_{2}$ and $C_{1} \subseteq C_{3} \subseteq \cdots \subseteq C_{2(m + 1)} \subseteq C_{2(m + 1) + 1}$. We will show that $C_{k} \cap C_{k + 1} = \emptyset$ for $0 \leq k \leq 2(m + 1)$. By the construction, $x_{1}, ...,  x_{t} \notin C_{2(m + 1)}$. Thus, $C_{2(m + 1)} \cap C_{2(m + 1) + 1} = \emptyset$. Hence, assume that $k < 2(m + 1)$. Observe that $C_{2i} = A \setminus \cup^{i}_{j = 1}A_{j} = A \setminus \{x_{1}, ..., x_{ip}\}, C_{2i + 1} = B \setminus \cup^{m - i + 1}_{j = 1}B_{j} = B \setminus \{x_{(i - 1)p + q + 1}, ..., x_{t}\}$ and $C_{2(i + 1)} = A \setminus \{x_{1}, ..., x_{(i +1)p}\}$. Since $(i - 1)p + q + 1 \leq ip$ and $(i - 1)p + q + 1 \leq (i + 1)p$, it follows that $C_{2i} \cap C_{2i + 1} = \emptyset$ and $C_{2i + 1} \cap C_{2(i + 1)} = \emptyset$.
\vskip 5 pt

\indent Now we are ready to construct the sets $S_{0}, ..., S_{2(m + 1) + 1}$. Firstly, let $S_{0} =  A$ and $S_{2(m + 1) + 1} = C_{2(m + 1) + 1} = B$. To construct $S_{1}$ and $S_{2}$, we distinguish $2$ subcases.
\vskip 5 pt

\noindent \textbf{Subcase 2.1 :} $r - t < p - q$. So,  $C_{1} = \emptyset$. We choose $r$ elements from $[n] \setminus S_{0}$ to add in $C_{1}$ in order to obtain $S_{1}$. Because $|[n] \setminus S_{0}| = l + r$, there are $\binom{l + r}{r} > \binom{l}{p}$ possible sets. Without loss of generality, we let
\vskip 5 pt

\indent $S_{1} = \{a_{1}, ..., a_{r}\}$.
\vskip 5 pt

\noindent Observe that $|C_{2}| = r - p$ and $|C_{3}| = r - mp = q + r - t$. We choose $p$ elements from $[n] \setminus (S_{1} \cup C_{2} \cup C_{3})$ to add in $C_{2}$ in order to obtain $S_{2}$. Recall that $C_{2} \cap C_{3} = \emptyset$, $S_{1} \cap C_{2} = \emptyset$ and $|S_{1} \cup C_{3}| \leq |S_{1}| + |C_{3}|$. Hence, we have that $|[n] \setminus (S_{1} \cup C_{2} \cup C_{3})| \geq 2r + l - (r + (r - p) + (q + r - t)) = l - (r - t - (p - q)) > l$, because $r - t < p - q$. Thus, there are $\binom{l - (r - t - (p - q))}{p} > \binom{l}{p}$ possible sets. Hence, we can choose a set $S$ among these sets which $S \cup C_{2} \notin \mathcal{C}$. We let
\vskip 5 pt

\indent $S_{2} = S \cup C_{2}$.
\vskip 5 pt

\noindent \textbf{Subcase 2.2 :} $r - t \geq p - q$. In this case, $|C_{1}| = r - (m + 1)p$. We choose $(m + 1)p$ elements from $[n] \setminus S_{0}$ to add in $C_{1}$ in order to obtain $S_{1}$. Clearly, there are $\binom{n - r}{(m + 1)p} = \binom{l + r}{(m + 1)p} > \binom{l}{p}$ possible sets. Hence, there exist $a_{1}, ..., a_{(m + 1)p} \in [n] \setminus S_{0}$ such that $\{a_{1}, ..., a_{(m + 1)p}, z_{p - q + 1}, ..., z_{r - t}\} \notin \mathcal{C}$. We let
\vskip 5 pt

\indent $S_{1} = \{a_{1}, ..., a_{(m + 1)p}, z_{p - q + 1}, ..., z_{r - t}\}$
\vskip 5 pt

\indent In order to construct $S_{2}$, we choose $p$ elements from $[n] \setminus (S_{1} \cup C_{2} \cup C_{3})$ to add in $C_{2}$. Since $C_{2} \subseteq S_{0}$ and $S_{1} \cap S_{0} = \emptyset$, it follows that $S_{1} \cap C_{2} = \emptyset$. We know that $C_{3} \cap C_{2} = \emptyset$. Moreover $|C_{3} \setminus S_{1}| \leq p$, because $|C_{3} \setminus C_{1}| = p$ and $C_{1} \subseteq S_{1}$. Hence, $|[n] \setminus (S_{1} \cup C_{2} \cup C_{3})| \geq |[n]| - (|S_{1}| + |C_{2}| + |C_{3} \setminus S_{1}|) \geq 2r + l - (r + (r - p) + p) = l$. So, there are at least $\binom{l}{p}$ possible sets. If there exists a set $R \in \binom{[n] \setminus (S_{1} \cup C_{2} \cup C_{3})}{p}$ such that $R \cup C_{2} \notin \mathcal{C}$, then we let
\vskip 5 pt

\indent $S_{2} = R \cup C_{2}$.
\vskip 5 pt

\indent Hence, we may assume that $\mathcal{C} = \{R \cup C_{2}:$ for all $R \in \binom{[n] \setminus (S_{1} \cup C_{2} \cup C_{3})}{p}\}$. Since $|S_{0} \setminus C_{2}| = p < l \leq |[n] \setminus (S_{1} \cup C_{2} \cup C_{3})|$, we can choose $w_{1}, ..., w_{p}$ from $[n] \setminus (S_{1} \cup C_{2} \cup C_{3})$ such that $w_{1} \notin S_{0} \setminus C_{2}$. By the construction of $S_{1}$, $S_{1} \setminus C_{1} = \{a_{1}, ..., a_{(m + 1)p}\}$ and $|S_{1} \setminus C_{1}| = (m + 1)p$. Hence, when $m \geq 1$, we have $|C_{3} \setminus C_{1}| = p$. Thus, $|S_{1} \setminus C_{1}| > |C_{3} \setminus C_{1}|$. This implies that at least one element in $\{a_{1}, ..., a_{(m + 1)p}\}$ is not in $C_{3} \setminus C_{1}$. Without loss of generality, we let $a_{1} \notin C_{3} \setminus C_{1}$. We have that $\{a_{1}, w_{2}, ..., w_{p}\} \cup C_{2} \notin \mathcal{C}$. Thus, we can let
\vskip 5 pt

\indent $S_{2} = \{a_{1}, w_{2}, ..., w_{p}\} \cup C_{2}$ and change $S_{1}$ to be
\vskip 5 pt

\indent $S_{1} = \{w_{1}, a_{2}, ..., a_{(m + 1)p}, z_{p - q + 1}, ..., z_{r - t}\}$.
\vskip 5 pt

\noindent When $m = 0$, we have $t = q$ and $a_{(m + 1)p} = a_{p}$. Thus, $B = C_{3} = S_{3}$, moreover, $S_{1} = \{a_{1}, ..., a_{p},$ $ z_{p - q + 1}, ..., z_{r - t}\}$. If $a_{1}, ..., a_{p} \in C_{3} \setminus C_{1}$, then there exists $a_{i}$ such that $a_{i} \in S_{0}$ contradicting the choice of $S_{1}$. Thus, there exists $a_{i} \notin C_{3} \setminus C_{1}$ and we can construct $S_{1}, S_{2}$ by similar arguments as the case when $m \geq 1$. This proves Case 2.
\vskip 5 pt

\indent We now construct $S_{3}, ..., S_{2(m + 1)}$. First we construct $S_{2i + 1}$ for $1 \leq i \leq m$. By the construction of $C_{2i + 1}$, we have $|C_{2i + 1}| = r - p(m - i + 1)$. Thus, we choose $p(m - i + 1)$ elements from $[n] \setminus (S_{2i} \cup C_{2i + 1})$ to add in $C_{2i + 1}$ in order to obtain $S_{2i + 1}$. Clearly, $|[n] \setminus (S_{2i} \cup C_{2i + 1})| = 2r + l - (r + (r - p(m - i + 1))) = l + p(m - i + 1)$. Hence, there are $\binom{l + p(m - i + 1)}{p(m - i + 1)} > \binom{l}{p}$ possible sets. So, we can choose a set $T$ among these sets which $T \cup C_{2i + 1} \notin\mathcal{C}$. We let
\vskip 5 pt

\indent $S_{2i + 1} = T \cup C_{2i + 1}$.
\vskip 5 pt

\noindent Next we will construct $S_{2i}$ where $2 \leq i \leq m + 1$. By the above construction, when we construct $S_{2i - 1}$, if all the $p(m - (i - 1) + 1)$ elements that we choose from $[n] \setminus (S_{2(i - 1)} \cup C_{2i - 1})$ to add in $C_{2i - 1}$ are not in $C_{2i + 1}$, then $|C_{2i + 1} \setminus S_{2i - 1}| = p$. But if some of which are in, then $|C_{2i + 1} \setminus S_{2i - 1}| < p$. In both case, $|C_{2i + 1} \setminus S_{2i - 1}| \leq p$. In order to obtain $S_{2i}$, we choose $ip$ elements from $[n] \setminus (S_{2i - 1} \cup C_{2i} \cup C_{2i + 1})$ to add in $C_{2i}$. Recall that $C_{2i} \cap C_{2i + 1} = \emptyset$.
\vskip 5 pt

\noindent We consider the case when $r - t \geq p - q$. Thus, $|C_{2i}| = r - ip$ for all $2 \leq i \leq m + 1$ implying that $|S_{2i - 1} \cup C_{2i} \cup C_{2i + 1}| \leq |C_{2i}| + |S_{2i - 1}| + |C_{2i + 1} \setminus S_{2i - 1}| \leq (r - ip) + r + p$. Hence, $|[n] \setminus (S_{2i - 1} \cup C_{2i} \cup C_{2i + 1})| \geq 2r + l - ((r - ip) + r + p) = l + (i - 1)p$. There are $\binom{l + (i - 1)p}{ip} > \binom{l}{p}$ possible sets. So, we can choose a set $W$ among these sets which $W \cup C_{2i} \notin \mathcal{C}$. We let
\vskip 5 pt

\indent $S_{2i} = W \cup C_{2i}$.
\vskip 5 pt

\noindent Hence, we consider the case when $r - t < p - q$. For $2 \leq i \leq m$, we can construct $S_{2i}$ by the same arguments as when $r - t \geq p - q$. Then, we may assume that $i = m + 1$. So $C_{2m + 2} = \emptyset$. Therefore, we choose $r$ elements from $[n] \setminus (S_{2m + 1} \cup C_{2m + 2} \cup C_{2m + 3})$. We see that $|(S_{2m + 1} \cup C_{2m + 2} \cup C_{2m + 3})| \leq |C_{2m + 3}| + |S_{2m + 1}| + |C_{2m + 2} \setminus S_{2m + 1}| \leq r + p$ which implies that $|[n] \setminus (S_{2m + 1} \cup C_{2m + 2} \cup C_{2m + 3})|$ $ \geq 2r + l - (r + p) = r + l - p > l$. There are at least $\binom{r + l - p}{r} > \binom{l}{p}$ possible sets. So, there is a set $U$ among these sets which is not in $\mathcal{C}$. We let
\vskip 5 pt

\indent $S_{2i} = U$.
\vskip 5 pt

\noindent By the choice of $S_{k}$, we see that $S_{k} \cap S_{k + 1} = \emptyset$ for all $0 \leq k \leq 2(m + 1)$. This proves the claim.
%\smallqed
\qed

\indent Now, we may let $A \in \mathcal{A}$ and $B \in \mathcal{B}$. By Claim \ref{absorb}, there exist $S_{0}, S_{1}, ..., S_{f}(f \geq 1)$ such that $A = S_{0}, B = S_{f}$, and for $k = 0, ..., f - 1, S_{k} \cap S_{k + 1} = \emptyset$ and $S_{k}, S_{k + 1} \notin \mathcal{C}$.
\vskip 5 pt

\indent We will show that $S_{k} \in \mathcal{A}$ for all $0 \leq k \leq f$. Clearly, $S_{0} = A \in \mathcal{A}$. Let $q = max\{i : S_{0}, ..., S_{i} \in \mathcal{A}\}$. We assume to the contrary that $q < f$. Remind that $\mathcal{A} \cup \mathcal{B} = \binom{[n]}{r} \setminus \mathcal{C}$ and $S_{q}, S_{q + 1} \notin \mathcal{C}$. Because $S_{q} \cap S_{q + 1} = \emptyset$ and the families $\mathcal{A}, \mathcal{B}$ are cross intersecting, it follows that $S_{q + 1} \in \mathcal{A}$ contradicts the maximality of $q$. Thus, $S_{k} \in \mathcal{A}$ for all $0 \leq k \leq f$, in particular, $B = S_{f} \in \mathcal{A} \cap \mathcal{B}$. This contradicts the disjoint property of $\mathcal{A}$ and $\mathcal{B}$. Hence, $|\mathcal{A}| + |\mathcal{B}| < \binom{n}{r} - \binom{l}{p}$. This completes the proof.
%\qed
\qed

\end{document}